\documentclass[11pt]{amsart}
\usepackage[sorted-cites]{amsrefs}
\usepackage{verbatim}
\usepackage{color}
\usepackage{amscd}

\newcommand{\overbar}{\overline}

\newfont{\fnt}{cmr10 scaled 550}

\newtheorem{theorem}{Theorem}

\newtheorem{lemma}{Lemma}
\newtheorem{cor}{Corollary}
\newtheorem{prop}{Proposition}

\newtheorem{definition}{Definition}

\theoremstyle{remark}
\newtheorem{remark}{Remark}

\font\strange=msbm10

\renewcommand{\epsilon}{\varepsilon}

\renewcommand{\Sigma}{\varSigma}
\newcommand{\R}{{{\mathchoice  {\hbox{$\textstyle{\text{\strange R}}$}}
{\hbox{$\textstyle{\text{\strange R}}$}}
{\hbox{$\scriptstyle  N\kern-0.3em  R$}}
{\hbox{$\scriptscriptstyle  R\kern-0.2em  R$}}}}}

\renewcommand{\phi}{\varphi}
\usepackage{amsmath,amsthm,amscd}

\begin{document}

\title[]{Uniqueness  of grim hyperplanes for mean curvature flows}



\thanks{The second author was partially supported by CNPq and Faperj of Brazil.}

\address{Instituto de Matem\'atica e Estat\'\i stica, Universidade Federal Fluminense,
Niter\'oi, RJ 24020, Brazil}
\author[Ditter Tasayco ]{Ditter Tasayco}
\email{ditter.y.t@gmail.com}
\author[Detang Zhou]{Detang Zhou}
\email{zhou@impa.br}

\newcommand{\M}{\mathcal M}

\begin{abstract}In this paper we show that an immersed nontrivial translating soliton for mean curvature flow  in $\mathbb{R}^{n+1}$($n=2,3)$  is a grim hyperplane if and only if it is mean convex and has  weighted total extrinsic curvature of at most quadratic growth. For  an embedded  translating soliton $\Sigma$ with nonnegative scalar curvature, we prove that if the mean curvature of $\Sigma$ does not change signs on each end, then $\Sigma$ must have positive scalar curvature unless it is either a hyperplane or a grim hyperplane.
\end{abstract}

\maketitle
\section{Introduction}

A mean curvature flow (MCF) in $\mathbb{R}^{n + 1}$ is  the negative gradient flow of the volume functional, which can be analyzed from the perspective of partial differential equations as shown by Huisken in \cite{huisken}. 
MCF is smooth in a short time and singularities must happen over a longer time. According to the rate of blow-up of the second fundamental form $A\left(t,p\right)$ of the hypersurface $\Sigma_t$, this finite time singularity $T$ is called type-I, if there exists a constant $C_0$ such that
\begin{equation*}
\displaystyle \sup_{p \in \Sigma_t}\left|A\left(t,p\right)\right|^2 \leq \frac{C_0}{\left(T - t\right)}
\end{equation*} 
for all $t < T$. Otherwise this finite time singularity is called type-II. 

We will deal with translating solitons which are important in study of type-II singularities. 

A complete connected isometrically immersed hypersurface $(\Sigma,\Phi)$ in $\mathbb{R}^{n+1}$ is called a {\it translating soliton}  if its mean curvature vector satisfies
\begin{equation*}
 \vec{H} =w^{\bot},
\end{equation*}
where $w \in {\mathbb{R}}^{n + 1}$ is a unitary vector and ${w}^{\bot}$ stands for the orthogonal projection of $w$ onto the normal bundle of $\Phi$. Let $\nu$ denote the unit normal along $\Phi$, then it is equivalent to
\begin{equation*}
\displaystyle H = -\left\langle \nu,w\right\rangle. \label{2I}
\end{equation*}

In particular, considering $f : \mathbb{R}^{n + 1} \longrightarrow \mathbb{R}$ defined by $f\left(x\right) = - \left\langle x,w\right\rangle$, then $\overbar{\nabla}f = -w$ and $H = \left\langle \overbar{\nabla}f,\nu\right\rangle$, therefore by definition translating solitons are $f$-minimal hypersurfaces. Since MCF is invariant under isometries, without loss of generality we may suppose that $w = \left(0,\ldots,0,1\right)$, then the function $f$ is defined by $f\left(x\right) = -x_{n + 1}$ and the $L_f$-stability operador of $\Sigma$ is given by
\begin{equation}
\displaystyle L_f =\Delta_f + \left|A\right|^2 \label{3I}  
\end{equation}

There are some examples of translating solitons: vertical hyperplanes, grim hyperplanes, translating bowl solitons and translating catenoids. In this article we will give a characterization of grim hyperplanes in dimensions $2$ and $3$. 

Recall that a {\it grim hyperplane} in $\mathbb{R}^{n+1}$ is a hypersurface $\mathcal{G}$ of $\mathbb{R}^{n+1}$ which can be represented parametrically via the embedding $\mathbf{\Phi} : \left(-\frac{\pi}{2},\frac{\pi}{2}\right)\times\mathbb{R}^{n - 1} \longrightarrow \mathbb{R}^{n + 1}$ defined by
\begin{eqnarray*}
\displaystyle 
\mathbf{\Phi}\left(t,y_1,\ldots,y_{n - 1}\right) = \left(t,y_1,\ldots,y_{n - 1},-\ln\left(\cos t\right)\right).
\end{eqnarray*}

The grim hyperplane $\mathcal{G}$ satisfies the translating soliton equation with $\displaystyle w = \left(0,\ldots,0,1\right)$ i.e. it is $f$-minimal for $\displaystyle f\left(x_1,\ldots,x_{n + 1}\right) = -x_{n + 1}$. Also it has positive mean curvature. When $n=2$ or $3$,  there exists a constant $C > 0$ such that 
\begin{equation}
\displaystyle \int_{B_R}{\left|A\right|}^2e^{-f} \ \leq \ C R^2 \label{2.29E} 
\end{equation}
for all $R$ sufficiently large. The aim of this article is to prove that indeed the grim hyperplanes are the only ones with these properties when $n=2, 3$.

\begin{theorem}\label{thm1} 
Let $\Phi : \Sigma^n \longrightarrow \mathbb{R}^{n + 1}$ be a translating soliton, with $n=2$ or $3$,  which is not a hyperplane. Then $\Sigma$ is a grim hyperplane if and only if $H=-\langle w,\nu\rangle \geq 0$ and there exists $C > 0$ such that
\begin{equation}
\displaystyle \int_{B_R}{\left|A\right|}^2e^{-f} \ \leq CR^2, \label{2.12B} 
\end{equation}
for all $R$ sufficiently large, where $B_R$ is the geodesic ball of radius $R$ and $f(x)=-\langle x,w\rangle$.
\end{theorem}

The expression (\ref{2.12B}) is not satisfied for $n \geq 4$ (see Proposition \ref{prop1}), thus Theorem \ref{thm1} is sharp in this sense.

It has been known that if $H\ge 0$ on a translating soliton $\Sigma$, then either $H\equiv 0$ on $\Sigma$ and $\Sigma $ is a hyperplane,  or $H>0 $ everywhere on $\Sigma$. 
Note that both hyperplane and grim hyperplane has vanishing scalar curvature. In \cite {MR3412395}, Mart{\'{\i}}n-Savas-Halilaj-Smoczyk proved that flat hyperplane and grim hyperplane are the only translating soliton with vanishing scalar curvature. It would be interesting to ask if the following is true.

{\bf Problem:}
 Let    $\Sigma$  be a translating soliton with nonnegative scalar curvature $S$.  Is it true that either $S\equiv 0$ on $\Sigma$ and   $ \Sigma $ is a hyperplane or grim hyperplane,  or $S>0 $ everywhere on $\Sigma$?

This problem is related to a result proved by Huang-Wu in \cite{MR3128984}.
 Let   $M$ be a closed embedded $n$-dimensional hypersurface in $\mathbb{R}^{n+1}$ with nonnegative scalar curvature. Let $M_t$ be a solution to the mean curvature flow with initial hypersurface $M$. Then the scalar curvature of $M_t$ is strictly positive for all $t>0$. 

 For complete embedded  translating solitons, we have 

\begin{theorem}\label{thm2}
Let $(\Sigma^{n}, g)$ be a embedded translating soliton with nonnegative scalar curvature $S$. Assume $H$ does not change signs on each end.
Then either $\Sigma$ is  a hyperplane or a grim hypersurface; or $\Sigma $ has positive scalar cuvature.
\end{theorem}

\section{Total weighted extrinsic curvature}
In this section we will give the asymptotic properties of the total weighted extrinsic curvatures of grim hyperplanes.
We have $$\partial_t = \sec\left(t\right)\left(\cos t, 0,\cdots, 0, \sin t\right).$$
We choose the unit normal $\nu$ to $\mathcal{G}$ to be
$\displaystyle \nu = \left(\sin t,0,\cdots,0,-\cos t\right).$
A little computation shows that $\displaystyle {\overbar{\nabla}}_{\partial_t} \nu = \left(\cos t\right) \partial_t$ and $\displaystyle {\overbar{\nabla}}_{\partial_{y_i}} \nu = 0$ \ $\displaystyle \left(1 \leq i \leq n - 1\right)$. 

Then the principal curvatures are $\lambda_1 = \cos t$, $\lambda_2 = \ldots = \lambda_n = 0$, thus on the coordinates $t$, $y_1$,$\ldots$, $y_{n - 1}$ the mean curvature only depends on $t$ and is given by $\displaystyle H\left(t\right) = \cos t$.
Since $t \in \left(-\frac{\pi}{2},\frac{\pi}{2}\right)$, we have the norm of the second fundamental form is given by
\begin{equation}
\displaystyle \left|A\right|\left(t\right)=\cos t=H\left(t\right).
\end{equation}
Now, consider the function $f : \mathbb{R}^{n + 1} \longrightarrow \mathbb{R}$ defined by $f\left(x\right) = -x_{n + 1}$, then
\begin{equation*}
\left\langle \overbar{\nabla}f,\nu\right\rangle =\cos t=H.
\end{equation*}

\begin{prop} \label{prop1}
The Grim Hyperplane $\mathcal{G}$ in $\mathbb{R}^{n + 1}$ satisfies
\begin{equation*}
\displaystyle \lim_{R \longrightarrow +\infty}\frac{1}{R^{n - 1}}\int_{B_R} \left|A\right|^2 e^{-f} = \left|B^{n - 1}\left(1\right)\right|\pi,
\end{equation*}
where $B_R$ is the geodesic ball with center at $0$ and radius $R$ and $\displaystyle B^{n - 1}\left(1\right)$ is the open ball in $\mathbb{R}^{n - 1}$ of radius $1$ and center at the origin.
\end{prop} 
{\bf Proof of Proposition \ref{prop1}.}
Observe that $f$ and the metric on $\mathcal{G}$ in the coordinates $t$, $y_1,\ldots,y_{n - 1}$ are given by
\begin{equation*}
\displaystyle f\left(t\right) = \ln\left(\cos t\right)
\end{equation*}
and
\begin{equation*}  g = \sec^2\left(t\right)dt^2 + dy^2_1 + \ldots + dy^2_{n - 1}.
\end{equation*}
Thus 
\begin{eqnarray*}
\displaystyle r =\int_{0}^{t} \sec\left(\xi\right)d\xi =-\ln\left(\tan\left(\frac{1}{2}\left(\frac{\pi}{2} - t\right)\right)\right), 
\end{eqnarray*}
we have $t = \frac{\pi}{2} - \eta\left(r\right)$, where $\eta\left(r\right) = 2\arctan\left(e^{-r}\right)$. Then %
\begin{equation*}
\displaystyle g = dr^2 + dy^2_1 + \cdots +dy^2_{n - 1}.
\end{equation*}
Besides that $\left|A\right|$ and $f$ in the coordinates $r, y_1,\cdots,y_{n - 1}$ are given by
\[
\displaystyle \left|A\right|\left(r\right) = \sin\left(\eta\left(r\right)\right),
\]
and
\[  f\left(r\right) = \ln\left(\sin\left(\eta\left(r\right)\right)\right).
\]
Denoting by $\left\|.\right\|$ the standard norm of $\mathbb{R}^{n - 1}$, we have
\begin{eqnarray*}
\displaystyle B_R &=& \left\{\left(r,y\right) \in \mathbb{R}\times\mathbb{R}^{n - 1} : r^2 + {\left\|y\right\|}^2 \leq R^2\right\} \\
\displaystyle &=& \left\{\left(r,y\right) \in \mathbb{R}\times\mathbb{R}^{n - 1} : -\sqrt{R^2 - {\left\|y\right\|}^2} \leq r \leq \sqrt{R^2 - {\left\|y\right\|}^2}, \quad \quad \left\|y\right\| \leq R\right\}.
\end{eqnarray*}
Since $-\eta'(r)=\sin\left(\eta(r)\right)$ is an even function, then
\begin{eqnarray*}
\displaystyle \int_{B_R}\left|A\right|^2 e^{-f} &=& \int_{\left\{\left\|y\right\| \leq R\right\}}\left[\int_{-\sqrt{R^2 - {\left\|y\right\|}^2}}^{\sqrt{R^2 - {\left\|y\right\|}^2}}\sin\left(\eta\left(r\right)\right)dr\right]dy \nonumber \\
 &=& \int_{\left\{\left\|y\right\| \leq R\right\}}\left[\pi - 2\eta\left(\sqrt{R^2 - \left\|y\right\|^2}\right)\right]dy \nonumber \\ 
\displaystyle &=& \pi \int_{\left\{\left\|y\right\| \leq R\right\}} 1 dy - 2\int_{\left\{\left\|y\right\| \leq R\right\}}\eta\left(\sqrt{R^2 - \left\|y\right\|^2}\right)dy \nonumber \\
\displaystyle &=& \pi\left|B^{n - 1}(1)\right|R^{n - 1} - 2\int_{0}^{R}\left(\int_{\mathbb{S}^{n - 2}_{\rho}}\eta\left(\sqrt{R^2 - \rho^2}\right)dA\right)d\rho \nonumber \\
\displaystyle &=& \pi\left|B^{n - 1}(1)\right|R^{n - 1} - 2\mbox{area}\left(\mathbb{S}^{n - 2}\right)\int_{0}^{R}\eta\left(\sqrt{R^2 - \rho^2}\right)\rho^{n - 2}d\rho.
\end{eqnarray*}
where we have used the co-area formula. Now, letting $\rho = R\sin\theta$ and using the fact $\mbox{area}\left(\mathbb{S}^{n - 2}\right) = (n - 1)\left|B^{n - 1}(1)\right|$, we have
\begin{equation}
\displaystyle \frac{1}{R^{n - 1}}\int_{B_R} \left|A\right|^2 e^{-f} = \left|B^{n - 1}\left(1\right)\right|\left[\pi - 2\left(n - 1\right)F_{n - 1}\left(R\right)\right], \label{5E} 
\end{equation}
where
\begin{equation}
\displaystyle F_{n - 1}\left(R\right) = \int_{0}^{\pi/2}\eta\left(R\cos\theta\right)\sin^{n - 2}\theta \cos\theta d\theta. \label{6E}
\end{equation}
Observe that
\begin{equation*}
\displaystyle \lim_{R \longrightarrow +\infty}\eta\left(R\cos\theta\right)\sin^{n - 2}\theta\cos\theta = 0 \qquad \mbox{for all} \ \theta \in \left[0,\frac{\pi}{2}\right].
\end{equation*}
Fixing $R > 0$, we have $\displaystyle \left|\eta\left(R\cos\theta\right)\sin^{n - 2}\theta\cos\theta\right| \leq  \frac{\pi}{2}\sin^{n - 2}\theta\cos\theta$ for all $\theta \in \left[0,\pi/2\right]$. Besides that 
\begin{equation*}
\int_{0}^{\pi/2}\sin^{n - 2}\theta\cos\theta d\theta = 1/(n - 1).
\end{equation*}
Then $\displaystyle \lim_{R \longrightarrow +\infty}F_{n - 1}\left(R\right) = 0$, and hence by (\ref{5E}), we get
\begin{equation*}
\displaystyle \lim_{R \longrightarrow +\infty}\frac{1}{R^{n - 1}}\int_{B_R}\left|A\right|^2 e^{-f} = \left|B^{n - 1}\left(1\right)\right|\pi.
\end{equation*}

\section{Proofs of Theorem \ref{thm1} and Theorem \ref{thm2}}
We begin this section with the following lemma which is in a form more general than we need. The lemma may have its independent interest.
\begin{lemma} \label{lem1} 
Assume that on a complete weighted manifold $\displaystyle \left(M,\left\langle ,\right\rangle,e^{-f}d\textsc{\footnotesize vol}\right)$, the functions $u$, $v \in C^2\left(M\right)$, with $u > 0$ and $v\ge 0$ on $M$, satisfy
\begin{eqnarray}
\displaystyle {\Delta}_f u + q\left(x\right) u \le 0 \qquad \qquad \mbox{and} \qquad \qquad {\Delta}_f v + q\left(x\right) v \ge 0, \label{eq11} 
\end{eqnarray}
where $q\left(x\right) \in C^{0}\left(M\right)$. Suppose that there exists a positive  function $\kappa > 0$  on $\mathbb{R}^+$ satisfying  $\frac{t}{\kappa(t)}$ is nonincreasing and 
\[\int^{+\infty}\frac{t}{\kappa(t)}dt=+\infty,
\] such that
\begin{equation}
\displaystyle \int_{B_R}v^2 e^{-f} \ \leq \kappa(R) \label{eq12}	
\end{equation} 
for all $R$. Then there exists a constant $C$ such that $\displaystyle v = Cu$.
\end{lemma}
\begin{remark} \label{remark1}
Without loss of generality, we can assume $\kappa(t)\geq C(1+t^2)$. Some examples of $\kappa(t)$ are $Ct^2$, $Ct^2\log(1+t)$, $Ct^2\log(1+t)\log\log(3+t)$,$\cdots$.
\end{remark}

{\bf Proof of Lemma \ref{lem1}.}
Set $ w = \frac{v}{u}$, then $v = wu$, thus by (\ref{eq11}) we get
\begin{equation*}
\begin{split}
\Delta_f v=&w\Delta_f u + 2\left\langle \nabla w, \nabla u\right\rangle+u \Delta_f w\\
\le&-w(qu)+2\left\langle\nabla w, \nabla u\right\rangle+u \Delta_f w \\
=&-qv + 2\left\langle \nabla w, \nabla u\right\rangle + u \Delta_f w.
\end{split}
\end{equation*}
Then
\begin{equation}
\displaystyle \Delta_f w \geq -2\left\langle \nabla w,\nabla\left(\ln u\right)\right\rangle. \label{3}	
\end{equation}

On the other hand, let $\phi \in C^2_o\left(M\right)$, then by (\ref{3}), we have
\begin{eqnarray*}
\displaystyle \int_{M}\phi^2{\left|\nabla w\right|}^2e^{-f} &=& \int_{M}\left\langle \phi^2\nabla w,\nabla w\right\rangle e^{-f} \\
\displaystyle &=& \int_{M}\left\langle \nabla\left(\phi^2 w\right),\nabla w\right\rangle e^{-f} - 2\int_{M}\phi w\left\langle \nabla\phi,\nabla w\right\rangle e^{-f} \\
\displaystyle &=& -\int_{M}\phi^2 w\left(\Delta_f w\right)e^{-f} - 2\int_{M}\phi w\left\langle \nabla\phi,\nabla w\right\rangle e^{-f} \\
\displaystyle &\leq& 2\int_{M}\phi^2 w\left\langle \nabla w,\nabla\left(\ln u\right)\right\rangle e^{-f} - 2\int_{M}\phi w\left\langle \nabla\phi,\nabla w\right\rangle e^{-f} \\
\displaystyle &=& 2\int_{M}\left\langle \phi\nabla w, w\left(\phi\nabla\left(\ln u\right) - \nabla\phi\right)\right\rangle \\
\displaystyle &\leq& \frac{1}{2}\int_{M}\phi^2{\left|\nabla w\right|}^2e^{-f} + 2\int_{M}w^2{\left|\phi\nabla\left(\ln u\right) - \nabla\phi\right|}^2 e^{-f}.
\end{eqnarray*}

Then
\begin{equation}
\displaystyle \int_{M}{\phi}^2{\left|\nabla w\right|}^2 e^{-f} \leq 4\int_{M}w^2{\left|\phi\nabla\left(\ln u\right) - \nabla\phi\right|}^2 e^{-f} \qquad \forall \ \phi \in C^2_o\left(M\right). \label{15e}
\end{equation}

If $\displaystyle \psi \in C^{\infty}_o\left(M\right)$, then $\phi = \psi u \in C^2_o\left(M\right)$. Besides that, a little computation shows 
\begin{equation*}
\displaystyle \displaystyle \phi \nabla\left(\ln u\right) - \nabla \phi = - \left(\nabla \psi\right) u,
\end{equation*}

Thus, from (\ref{15e}), we have
\begin{eqnarray}
\displaystyle \int_{M}{\psi}^2 u^2{\left|\nabla w\right|}^2 e^{-f} &\leq& 4\int_{M}w^2\left|\nabla\psi\right|^2 u^2 e^{-f} \nonumber \\
\displaystyle &=& 4\int_{M}\left|\nabla\psi\right|^2 v^2 e^{-f} \qquad \forall \ \psi \in C^{\infty}_o\left(M\right). \label{2.20D}  
\end{eqnarray}
Define  functions $\beta$, $\xi $ on $[0,+\infty)$
as 
\[\beta(t):=\int_0^t\frac{\tau}{\kappa(\tau)}d\tau,
\]
and $\xi$ is the inverse function of $\beta$. From the hypothesis we know  $\beta'$  is nonincreasing and $\xi'$ is nondecreasing functions on $[0,+\infty)$.
Now, we now choose a cutoff  function 
\[
\psi_R(x) =
 \begin{cases}
 1, & \mathrm{ on } \ B_{\xi(R)} ;\\
2 - \frac{\beta(r(x))}{R},  & \mathrm{ on } \ B_{\xi(2R)} \setminus B_{\xi(R)};\\
0,& \mathrm{ on } \ M\setminus B_{\xi(2R)}.
\end{cases}
\]
where $r\left(x\right) = d\left(x,p\right)$, $p \in M$ is a fixed point and $B_R$ is the geodesic ball with radius $R$ and center $p$. We see that  $\displaystyle \left|\nabla\psi_R\right| = \frac{\beta'(r)}{R }=\frac{r}{R\kappa(r)}$.
Then, by (\ref{eq12}), we get
\begin{eqnarray*}
\displaystyle \int_{B_{\xi(R)}}u^2{\left|\nabla w\right|}^2 e^{-f} &=& \int_{B_{\xi(R)}}{\psi}^2_R u^2{\left|\nabla w\right|}^2 e^{-f} \\ 
\displaystyle &\leq& \int_{M}{\psi}^2_R u^2\left|\nabla w\right|^2 e^{-f} \\
\displaystyle &\leq& 4\int_{M} v^2 \left|\nabla \psi_R\right|^2 e^{-f} \\ 
\displaystyle &=& 4\int_{B_{\xi(2R)} \setminus B_{\xi(R)}} v^2 \left|\nabla \psi_R\right|^2 e^{-f} \\
\displaystyle &=& \frac4{R^2}\int^{\xi(2R)}_{\xi(R)}(\beta'(s))^2\int_{\partial B_s}v^2  e^{-f} dAds.
\end{eqnarray*}
Here we have used co-area formula. For convenience, we write $V(s)=\int_{ B_s}v^2  e^{-f} dV$. Therefore 
\[
V(s)=\int_0^s\int_{ \partial B_\tau}v^2  e^{-f} dAd\tau\le \kappa(s), 
\]
and 
\[
\begin{split}
\int_{B_{\xi(R)}}u^2{\left|\nabla w\right|}^2 e^{-f} \le& \frac4{R^2}\int^{\xi(2R)}_{\xi(R)}(\beta'(s))^2V'(s)ds\\
=&\frac{4}{R^2}\left[V(s)(\beta'(s))^2\left|^{\xi(2R)}_{\xi(R)}\right.-\int^{\xi(2R)}_{\xi(R)}2V(s)(\beta'(s))d\beta'(s)\right]\\
\le &\frac{4}{R^2}\left[V(s)(\beta'(s))^2\left|^{\xi(2R)}_{\xi(R)}\right.-2\int^{\xi(2R)}_{\xi(R)}sd\beta'(s)\right]\\
\le &\frac{4}{R^2}\left[V(s)(\beta'(s))^2\left|^{\xi(2R)}_{\xi(R)}\right.-2s\beta'(s)\left|^{\xi(2R)}_{\xi(R)}\right.+2\int^{\xi(2R)}_{\xi(R)}\beta'(s)ds\right]\\
\le &\frac{4}{R^2}\left[V(s)(\beta'(s))^2\left|^{\xi(2R)}_{\xi(R)}\right.-2s\beta'(s)\left|^{\xi(2R)}_{\xi(R)}\right.+\beta(s)\left|^{\xi(2R)}_{\xi(R)}\right.\right]\\
\le &\frac{4}{R^2}\left[V(s)(\beta'(s))^2\left|^{\xi(2R)}_{\xi(R)}\right.-2s\beta'(s)\left|^{\xi(2R)}_{\xi(R)}\right.+R\right]
\end{split}
\]
Since
\[V(s)(\beta'(s))^2=V(s)\beta'(s)\beta'(s)\le s\beta'(s),
\]
and $\beta'(s)=\frac{s}{\kappa(s)}$, thus Remark \ref{remark1} implies these terms are bounded,
hence when $R \longrightarrow +\infty$,  all the terms on the right hand side go to zero. So we get
\begin{equation*}
\displaystyle \int_{M} u^2 {\left|\nabla w\right|}^2 e^{-f} = 0.
\end{equation*}
Then $\nabla w \equiv 0$, thus there is a constant $C$ such that $w \equiv C$ and hence $v =C u$.
\qed

\begin{definition}
A two-sided translating soliton $\Sigma$ is said to be stable if
\begin{equation*}
\displaystyle \int_{\Sigma}\left[{\left|\nabla\varphi\right|}^2 - |A|^2{\varphi}^2\right]{e}^{-f}d\sigma \geq 0 \quad \mbox{for all} \ \varphi \in {C}^{\infty}_{o}\left(\Sigma\right).
\end{equation*}
\end{definition}
As a consequence of Lemma \ref{lem1}, we have the following:

\begin{cor} \label{corollary2.13}  
Let $\Phi : \Sigma^n \longrightarrow \mathbb{R}^{n + 1}$ be a stable translating soliton and let $\omega \in C^{2}\left(\Sigma\right)$ be a positive solution of the stability equation
\begin{equation}
\displaystyle \Delta_f \omega  + {\left|A\right|}^2 \omega =0. \label{2.17}  
\end{equation}
Moreover,  if $H\ge 0$ and  there exists a constant $C > 0$ such that
\begin{equation}
\displaystyle \int_{B_R}H^2 e^{-f} \ \leq \ CR^2 \qquad \mbox{for all} \ R \ \mbox{large enough}. \label{2.18} 
\end{equation} 
Then there exists a constant $\widetilde{C}$ such that $H = \widetilde{C}\omega$. In particular, if $H \not\equiv 0$, then $\widetilde{C} \in \mathbb{R} \setminus \left\{0\right\}$ and  $H>0$.

\end{cor}

Now, we include here a  result due to Li and Wang (\cite{pw}) which will be needed in the proof our main theorem.
\begin{lemma} \label{lem2}
Suppose $\Sigma$ is complete and there exists a nonnegative function $\phi : \Sigma \longrightarrow \mathbb{R}$, not identically zero, such that $(\Delta_f+q)\left(\phi\right) \leq 0$. Then $\Delta_f+q$ is stable.
\end{lemma}
{\bf Proof.}
Let $\Omega$ be a compact subdomain in $\Sigma$ and let $u$ be the first eigenfunction satisfying
\begin{eqnarray}
\displaystyle 
\left\{
\begin{array}{rcl}
\displaystyle (\Delta_f+q) u =& - \lambda_1(\Omega) u  & \mbox{in} \ \Omega, \\
[6pt]
\displaystyle u = &0 &\mbox{on} \ \partial\Omega.
\end{array}
\right. \label{1.2D}  
\end{eqnarray}
We may assume that $u \geq 0$ on $\Omega$. From regularity of $u$ and Hopf Lemma, we have
\begin{eqnarray*}
\displaystyle &\bullet& u > 0 \ \ \mathrm{in \ the \ interior \ of} \ \Omega. \qquad \qquad \qquad \qquad \\
\displaystyle &\bullet& \frac{\partial u}{\partial \nu} < 0 \ \mathrm{on} \ \partial \Omega, \ \mathrm{where} \ \nu \ \mathrm{is \ the \ outward \ unit \ normal \ of} \ \partial\Omega. 
\end{eqnarray*}
Thus, integration by parts on $u$ and $\phi$ and also the hypothesis, we have
\begin{eqnarray}
\displaystyle \int_{\Omega}u\left(\Delta_f \phi\right)e^{-f} - \int_{\Omega}\phi \left(\Delta_f u\right)e^{-f} &=& \int_{\partial\Omega} u\frac{\partial \phi}{\partial \nu}e^{-f} - \int_{\partial\Omega} \phi\frac{\partial u}{\partial \nu}e^{-f} \nonumber \\
\displaystyle &=& -\int_{\partial\Omega} \phi\frac{\partial u}{\partial \nu}e^{-f} \geq 0. \label{1.3D}  
\end{eqnarray}

From hypothesis and (\ref{1.2D}), we have
\begin{eqnarray}
\displaystyle 
\left\{
\begin{array}{rcl}
\displaystyle \Delta_f \phi + Q\phi &\leq& 0, \\
\displaystyle \Delta_f u + Qu &=& -\lambda_1\left(\Omega\right)u.
\end{array}
\right.
\label{1.4D}  
\end{eqnarray}

Since $u > 0$, multiplying the first inequality of (\ref{1.4D}) by $u$ and the second equation by $-\phi$, and finally both by $e^{-f}$, we have
\begin{equation}
\displaystyle u\left(\Delta_f\phi\right)e^{-f} - \phi\left(\Delta_f u\right)e^{-f} \ \leq \ \lambda_1\left(\Omega\right) \left(\phi u\right)e^{-f} \label{1.5D}  
\end{equation}

Since both $u > 0$ and $\phi \geq 0$ are not identically zero, then combining (\ref{1.5D}) with (\ref{1.3D}), we have $\lambda_1\left(\Omega\right) \geq 0$ for all compact subdomains of $\Sigma$, then $\lambda_1\left(f,Q\right) \geq 0$, therefore $\Delta_f+q$ is stable.
\qed

We are now ready to give the proof of Theorem \ref{thm1}.

{\bf Proof of Theorem \ref{thm1}}
 Since $\Phi : \Sigma^n \longrightarrow \mathbb{R}^{n + 1}$ is a translating soliton, then the mean curvature $H$ satisfies $\displaystyle \Delta_f H + \left|A\right|^2 H = 0$(see Proposition 3 in \cite{CMZ}). Since $H \geq 0$ and $\Sigma$ is a non-planar translating soliton , then $H$ is not identically zero, thus by Lemma \ref{lem2}, $\Sigma$ is stable and hence the weighted version of a result by Fischer-Colbrie and Schoen \cite{FS} guarantees there exists a non-constant positive $C^2$-function $\omega$ on $\Sigma$ such that
\begin{equation}
\displaystyle \Delta_f \omega + {\left|A\right|}^2 \omega =0.
\end{equation}

As $\displaystyle \frac{H^2}{n} \leq {\left|A\right|}^2$ and $\left|A\right|$ satisfies (\ref{2.12B}), then
\begin{equation}
\displaystyle \int_{B_R}H^2 e^{-f} \ \leq \ nC R^2.
\end{equation}

Then, by Corollary \ref{corollary2.13} and the condition that $H \geq 0$ and not identically zero, there is a constant $C_1 > 0$ such that
\begin{equation}
\displaystyle H =C_1 \omega. \label{2.27B}  
\end{equation}

In particular $H > 0$ everywhere on $\Sigma$. On the other hand,  the Simons equation( see  \cite{CMZ} or \cite{MR3412395})  implies that
\begin{equation}
\displaystyle \left|A\right|\left\{\Delta_f \left|A\right| + {\left|A\right|}^2\left|A\right|\right\} ={\left|\nabla A\right|}^2 - {\left|\nabla \left|A\right|\right|}^2 \ \geq \ 0. \label{2.28B}  
\end{equation}

Since $\left|A\right|$ satisfies (\ref{2.12B}), then by Lemma \ref{lem1}, $\exists \ C_2 \geq 0$ such that
\begin{equation}
\displaystyle \left|A\right| =C_2 \omega. \label{2.34U}  
\end{equation}

Besides that $\Sigma^{n}$ is not a hyperplane, then $\left|A\right|$ is not identically zero, thus $C_2 > 0$. Then by (\ref{2.27B}) and (\ref{2.34U}) we have ${\left|A\right|}^2H^{-2} = \mathrm{constant} > 0$. In particular this function attains its local maximum on $\Sigma$.  Theorem B in \cite{MR3412395} says that $\Sigma$  is a grim hyperplane if and only if the function  ${\left|A\right|}^2H^{-2}$attains a local maximum. Therefore $\Sigma$ is a grim hyperplane. 
\qed

We now prove Theorem \ref{thm2}.

{\bf Proof of Theorem \ref{thm2}.} To prove Theorem \ref{thm2}, we will need a result of Huang-Wu\cite{MR3128984}.  Denote by $M_+$ a connected component of $\{p\in M,  H\ge 0 \textrm{ at } p\}$ that contains a point of positive mean curvature. We say that the mean curvature $H$ changes signs through $\Gamma$ if $\Gamma$ is a connected component of $\partial M_+$ and $\Gamma$ intersects the boundary of a connected component of $M\backslash\partial M_+$.    Theorem 2 of Huang-Wu\cite{MR3128984} $S\ge 0$, says that  if $H$  changes sign along $\Gamma$ then $\Gamma$ is  unbounded set. Since we have assumed that $H$ does not changes signs at infitiy,   $H$ has a sign. Hence either

(1) $H\equiv 0$, or 

(2)$H\ge 0$ but does not vanish at least one point. 

In case (1), $\Sigma $ must be a hyperplane.

In case (2), if there is point $p\in \Sigma$, such that $S(p)=0$ then $|A|^2=H^2-S\le H^2$ and equality holds at $p$. Therefore the function $|A|^2H^2$ is well defined and attains its maximum  at $p$. By Theorem B in \cite{MR3412395} it must be a grim hyperplane.
\qed

\end{document}